\documentclass[10pt,reqno]{amsart}
\usepackage[margin=1.75in,letterpaper]{geometry}   
\usepackage{setspace}
\usepackage{graphicx}
\usepackage{tikz}
\usepackage{tikz-cd}

\usetikzlibrary{graphs, positioning}
\usepackage{comment}
\usepackage{amsmath, amssymb, amsfonts, enumerate}
\usepackage{amsthm} % Theorem environments
\usepackage{mathtools} % Extended math tools
\usepackage[colorlinks=true, pdfstartview=FitV, linkcolor=blue, citecolor=blue, urlcolor=blue]{hyperref}
\usepackage[capitalise]{cleveref}
\crefformat{equation}{(#2#1#3)}
\crefrangeformat{equation}{(#3#1#4--#5#2#6)}
\crefformat{enumi}{(#2#1#3)}
\crefrangeformat{enumi}{(#3#1#4--#5#2#6)}

\makeatother

\newcounter{intro}

\newtheorem{theorem}[subsection]{Theorem}
\newtheorem{proposition}[subsection]{Proposition}
\newtheorem{lemma}[subsection]{Lemma}
\newtheorem{corollary}[subsection]{Corollary}

\theoremstyle{definition}
\newtheorem{example}[subsection]{Example}

\newtheorem{definition}[subsection]{Definition}

\newtheorem{chunk}[subsection]{}

\newtheorem{obs}[subsection]{Observation}

\theoremstyle{remark}
\newtheorem{remark}[subsection]{Remark}

\numberwithin{equation}{section}

\newcommand{\Tor}{\mathrm{Tor}}

\newcommand{\reg}{\mathrm{reg}}

\newcommand{\kk}{\mathsf k}

\newcommand{\m}{\mathfrak m}

\numberwithin{equation}{section}

\title[Product of ideals and Golod threshold]{Golod thresholds of ideals  and products of ideals in dimension 3}

\author[D.~Limbu]{Dinesh Limbu}
\address{Department of Mathematics, University of South Carolina, 
Columbia, SC 29208}
\email{dlimbu@email.sc.edu}

\author[A.~Vraciu]{Adela Vraciu}
\address{Department of Mathematics, University of South Carolina, 
Columbia, SC 29208}
\email{vraciu@math.sc.edu}

\date{\today}

\begin{document}

\begin{abstract}
Let $I$ be a proper ideal in a standard graded polynomial ring $P$ with homogeneous maximal ideal $\m$. We show that there exists $k_0\in \mathbb{N}$ such that $\m^kI$ is Golod for all $k\geq k_0$. The smallest such $k_0$ is denoted $\mathrm{gt}_{\m}(I)$, and we call it the Golod threshold for $I$ over $\m$. We give an explicit formula for $\mathrm{gt}_{\m}(I)$ in the case when $I$ is a monomial complete intersection and bounds for $\mathrm{gt}_{\m}(I)$ for other classes of ideals. We also show that when $I$ is homogeneous, under a suitable degree condition, $\m^k I$ is Golod if and only if $k\geq \mathrm{gt}_{\m}(I)$, i.e., once Golodness is achieved, it cannot be lost. In the last section, we work in the setting of a regular local ring of dimension 3, and show that the product of any two proper ideals is Golod; this resolves a question posed by Dao and De Stefani \cite[Question 4.5]{dao2022monomial} and VandeBogert \cite[Question 4.11]{vandebogert2022products}.
\end{abstract} 



\maketitle

\section{Introduction}
Let $(P,\m,\kk)$ be a regular local ring of dimension $n$ over a field $\kk$ or a standard graded polynomial ring in $n$ variables, with homogeneous maximal ideal $\m$. Let $I$ be a proper ideal of $P$ and let $R=P/I$. The Poincar\'{e} series of $\kk $ over $R$ is $P_\kk^R(t)=\sum_{i\geq0}\beta_i^R(\kk)t^i$, where $\beta_{i}^R(\kk)=\dim_{\kk}[\Tor_i^R(\kk,\kk)]$, and the Poincar\'{e} series of $R$ over $P$ is $P_R^P(t)=\sum_{i\geq0}\beta_i^P(R)t^i$, where $\beta_{i}^P(R)=\dim_k[\Tor_i^P(R,\kk)]$. Serre established the following coefficient-wise inequality:
\[
P_\kk^R(t)
\leq
\frac{(1+t)^n}
{1-t\bigl(P_R^P(t)-1\bigr)}.
\]
When the above upper bound is attained, we say the ring $R$ is Golod and the ideal $I$ is Golod. The question of whether the product of two proper ideals is Golod has been studied in several papers: \cite{herzog1979golodideale}, \cite{herzog2013ordinary}, \cite{de2016products}, \cite{herzog2016homology}, \cite{maleki2019golod}, \cite{herzog2018koszul} and \cite{vandebogert2022products}. We briefly review some of the known results on this question. In the setting of a standard graded polynomial ring of characteristic $0$, Herzog and Huneke \cite{herzog2013ordinary} introduced strongly Golod ideals and showed that powers and symbolic powers of graded ideals are strongly Golod and the product of two strongly Golod ideals is Golod. Later, Herzog and Maleki \cite{herzog2018koszul} introduced $d$-Golod ideals in arbitrary characteristic and showed that the product of two $d$-Golod ideals is Golod. However, due to the work of De Stefani, we know that in more than three variables there are ideals whose product is not Golod. The following is the famous counterexample:\begin{example}\cite[Example 2.1]{de2016products}
    Let $P=\kk[x_1,x_2,x_3,x_4]$, $I=(x_1^2, x_2^2, x_3^2, x_4^2)$ and $\m=(x_1, x_2, x_3, x_4)$. Then $\m I$ is not Golod. 
\end{example} 
In the above example, we notice that $\m^kI$ is Golod for all $k\geq 2$. To capture this phenomenon, we define an invariant which we call the Golod threshold, as follows. 
\begin{definition}\label{golodthresholddef}
    Let $I$, $J$ be ideals in $P=\kk[x_1,\cdots,x_n]$.
    The {\em Golod threshold of $I$ with respect to $J$} is
    $$
    \mathrm{gt}_J(I):=\mathrm{min}\{k \in \mathbb N\, | \, J^{\ell } I \ \mathrm{is \ 
    Golod\ for \ all \ } \ell \ge k\},
    $$
    provided that such a $k$ exists.
\end{definition}
In Section 3 we show that the Golod threshold always exists: see \cref{existence of golod threshold}. We will then focus on the case $J=\m$ and establish upper bounds for all homogeneous ideals $I$ [\cref{another bound}, \cref{upper bound}] and a lower bound for homogeneous Gorenstein ideals $I$ [\cref{lower bound for gorenstein}], based on the graded shifts in the minimal free resolution of the ideal over the polynomial ring. Using these bounds, we compute the Golod threshold with respect to $\m$ for certain Gorenstein ideals; see Example \cref{gorenstein example}. At the end of this section, under a suitable condition on the degrees of the minimal generators of a homogeneous ideal $I$, we show that if $I$ is Golod then $\m I$ is Golod; see \cref{equigenratedgolodThreshold}. As a consequence, under suitable degree conditions, $\m^k I$ is Golod if and only if $k\geq \mathrm{gt}_{\m}(I)$.  

In Section 4, we develop a combinatorial criterion for the existence of a monomial basis for the Koszul homology; see \cref{shellmonbasis}. We apply this criterion to show that ideals of the form $\m^k (x_1^{a_1},\cdots,x_n^{a_n})$ admit monomial bases for their Koszul homology; see \cref{monomialbasisforCI}. 

In Section 5, we establish an exact formula for the Golod threshold with respect to $\m$ for the monomial complete intersection $(x_1^{a_1},\ldots,x_n^{a_n})$; see \cref{preciseCIbound}.

In the last section, we focus on the case of dimension three, where several results are already known. Dao and De Stefani \cite{dao2022monomial} showed that in a standard graded three-variable polynomial ring, the product of any two proper monomial ideals is Golod. More recently, VandeBogert \cite{vandebogert2022products} showed that in a regular local ring of dimension three, the product of the maximal ideal with any proper ideal is Golod. The question of whether the product of two proper ideals in a regular local ring of dimension three is Golod had remained open. In \cref{maintheorem1} we prove that in a regular local ring of dimension three, the product of any two proper ideals is Golod. This resolves the questions of Dao and De Stefani\cite[Question 4.5]{dao2022monomial} and VandeBogert\cite[Question 4.11]{vandebogert2022products} in the positive direction.

\section{Preliminaries}
Let $P=\kk[x_1,\ldots,x_n]$ be a standard multigraded polynomial ring over a field $\kk$ and $\m=(x_1,\ldots,x_n)$ be the maximal homogeneous ideal in $P$. Let $I$ be a proper ideal of $P$ and let $R=P/I$.
% The Koszul complex is one of the fundamental constructions in commutative algebra and has numerous applications throughout the subject. For more details, see \cite{avramov1998infinite}.
The Koszul complex on $x_1,\cdots,x_n$ over $R$ is \[
\mathbf{K}^{R}(\mathbf{x})=R\langle dx_1,\cdots,dx_n|\;\partial(dx_i)=x_i\rangle
\]
where $dx_i$ are exterior variables of degree 1. For $F=\{i_1<\cdots<i_r\}\subseteq [n]$, we will denote $dx_{i_1}\cdots dx_{i_r}$ by $dx_F$. For simplicity, we will use the notation $\mathbf {K}^{R}:=\mathbf{K}^{R}(\mathbf{x})$ and we will refer to this as the Koszul complex of $R$. For a regular local ring we will use the same notation as above, but instead of variables we will use the minimal generators of the maximal ideal to define the Koszul complex. For an element $z\in \mathbf{K}_i^R$, $|z|:=i$ denotes the homological degree of $z$.

\begin{chunk}
    If $I$ is a monomial ideal, the multigrading of $P$ naturally induces a multigrading on $R$, which in turn induces a multigrading on the Koszul complex of $R$. For $\alpha\in \mathbb{N}^n$, the $\alpha$-multigraded piece of $\mathbf{K}^{R}(\mathbf{x})$, would be\[
    [\mathbf{K}^R(\mathbf{x})]_{\alpha}=\bigoplus_{\substack{F\subseteq\mathrm{Supp}(x^{\alpha})\\ x^{\alpha-1_F}\notin I}}\kk x^{\alpha-1_F}dx_F
 \] 
 with the differential
 \[
\partial(x^{\alpha-1_F}dx_F)=\sum_{\substack{i\in F\\ x^{\alpha-1_{F\setminus \{i\} }}\notin I }} (-1)^{\mathrm{pos}_F(i)-1}x^{\alpha-1_{F\setminus \{i\} }}dx_{F\setminus \{i\} }
\]
where $\mathrm{pos}_F(i)$ is the position of $i$ in $F$ in the natural order and for $G\subseteq [n]$, $1_{G}=(g_1,\cdots,g_n)$ where $g_i=1$ if $i\in G$ otherwise $g_i=0$. 
\end{chunk}
\begin{chunk}
     Here we recall some terminology related to simplicial chain complexes. For more details, we refer to \cite[Chapter 1]{miller2005combinatorial}. A collection $\Delta\subseteq 2^{[n]}$ is called a simplical complex on vertex set $[n]$ if $F'\subseteq F$ and $F \in \Delta$ implies $F'\in \Delta$. The elements of $\Delta$ are called faces, and the maximal faces are called facets. The dimension of a face $F$ is $\dim F=|F|-1$. 
      Let \[
    \mathbb{F}_i=\{F\in \Delta:\;|F|=i+1\},\quad \mathbb{F}_{-1}=\{\emptyset\} \]
    and let $\kk^{\mathbb{F}_i}$ be the $\kk$-vector space with basis $\mathbb{F}_i$. The reduced chain complex of $\Delta$ over $\kk$ is \[
        \widetilde{\mathbf{C}}_{\cdot}(\Delta;\kk):\quad
0 \longrightarrow \mathbf{\kk}^{\mathbb{F}_{n-1}}
\xrightarrow{}
\cdots
\xrightarrow{}
\mathbf{\kk}^{\mathbb{F}_{0}}
\xrightarrow{}
\mathbf{\kk}^{\mathbb{F}_{-1}}
\longrightarrow 0.
\]
where for $F\in \mathbb{F}_i$ we define the differential by $
\partial_i(F)=\sum_{j\in F}(-1)^{\mathrm{pos}_{F}(j)-1}(F\setminus \{j\})
$. The $i^{th}$ reduced homology of $\Delta$ over $\kk$ is defined as $
\widetilde{H}_i(\Delta;\kk)=\ker(\partial_i)/\mathrm{im}(\partial_{i+1}).$ We will denote the reduced cochain complex of $\Delta$ over $\kk$ by $\widetilde{\mathbb C}^{\cdot}(\Delta;\kk)$, which is defind as $ \widetilde{\mathbb C}^{\cdot}(\Delta;\kk)
:=\mathrm{Hom}_{\kk}(\widetilde{{C}}_{\cdot}(\Delta;\kk),\kk)$ and its $i$th reduced cohomology of $\Delta$ over $\kk$ is defined as $\widetilde{H}^{i}(\Delta;\kk)=\ker(\partial^{i})/\mathrm{im}(\partial^{i-1})$. For a face $F\in \widetilde{\mathbf{C}}_{\cdot}(\Delta;\kk)$, we denote its dual in $ \widetilde{\mathbf{C}}^{\cdot}(\Delta;\kk)$ by $F^*$. 
     
\begin{comment}
\begin{definition}[Shellable Simplicial Complex]

A d-dimensional simplicial complex is called 'pure' if its maximal simplices all have dimension 'd'. Let $\Delta$ be a finite or countably infinite simplicial complex. An ordering $C_1,C_2,\ldots$ of the maximal simplices of $\Delta$ is a 'shelling' if, for all $k=2,3,\ldots$, the complex :\[
B_k:=\Big(\bigcup_{i=1}^{k-1}C_i\Big)\cap C_k
\]
is pure and of dimension one smaller than $\dim C_k$.
\end{definition}

\begin{definition}[Shifted Simplicial Complex]
A simplicial complex on $n$ vertices is shifted if there exists a labeling of the nodes by $1$ through $n$ such that for any face $\{v_1, v_2, \cdots, v_d\}$, replacing any $v_i$ by a node
with a smaller label results in a collection which is also a face. 
\end{definition}
\end{comment}
\end{chunk}

\begin{chunk} We will now recall an equivalent characterization of Golodness. For more details on Golod rings, one can check \cite{avramov1998infinite}.
 Let $I$ be a proper or proper homogeneous ideal in a regular local ring or a standard graded ring $P$ with maximal or maximal homogeneous ideal $\m$, and let $R=P/I$. Let $\mathfrak B$ be a $k$-basis for $H_{\geq 1}(\mathbf{K}^{R}(\mathbf{x}))$. A trivial Massey operation on the Koszul complex of $R$ is a function \[
     \mu:\bigsqcup _{q\geq 1} \mathfrak B^q\xrightarrow[]{}\mathbf{K}^{R}(\mathbf{x})  \ \ \text{such that}
     \]
\begin{itemize}
         \item[(i)] For $h\in \mathfrak B$, $\mu(h)$ is a cycle representing the class of $h$.
         \item[(ii)] For $q\geq 2$ and $h_1,\cdots,h_q\in \mathfrak B$, we can choose $\mu(h_1,\cdots,h_q)$ such that,
         \begin{equation}\label{massey}
         \partial \mu(h_1,\cdots,h_q)=\sum_{j=1}^{q-1}\overline{\mu(h_1,\cdots,h_{j})}\mu(h_{j+1},\cdots,h_q)
         \end{equation} 
         where $\overline{a}=(-1)^{|a|+1}a$.
         Note that $|\mu(h_1, \ldots, h_q)|=\sum_{j=1}^q |h_i|+q-1$.
         \item[(iii)] $\mathrm{im}(\mu)\subseteq \m \mathbf{K}^R(\mathbf{x})$
     \end{itemize}
 The ideal $I$ is Golod if and only if $\mathbf{K}^{P/I}(\mathbf{x})$ admits a trivial Massey operation.
 \begin{remark}\label{notGolodcondition}
 If $\mathbf{K}^R({\bf x})$ admits a trivial Massey operation, then the product of any two elements in $H_{\ge 1}({\bf K}^R({\bf x}))$ is zero.
 \end{remark}  
 \begin{remark}\label{zeroproduct} Suppose that we can choose a basis $\mathfrak B$ for $H_{\geq 1}(\mathbf{K}^{R}(\mathbf{x}))$, such that for every $h\in \mathfrak B$, there is a cycle $z_h$, with $h=[z_h]$ and \[z_hz_{h'}=0\quad\text{for all }h,h'\in \mathfrak B.
 \]Then we may define $\mu(h)=z_h$ and $\mu(h_1,\cdots,h_p)=0$ for $p\geq 2$. One can check that this is indeed a trivial Massey operation. The following is a slightly more general version of this observation.
 \end{remark}
\end{chunk}

\begin{lemma}\label{trivialmassey1}
Let $p:=\mathrm{pd}_P(P/I)$. Assume that $\mathfrak B$ is a basis for $H_{\ge 1}(\mathbf{K}^{R}(\mathbf{x}))$, and for every $h\in \mathfrak B$ there is a cycle $z_h$ with $[z_h]=h$ such that $z_hz_{h'}=0$ for all $h, h'\in \mathfrak B$ with $|h|+|h'|\le p$.
 Then $I$ is Golod.
\end{lemma}
\begin{proof}
 We can define $\mu(h)=z_h$ and $\mu(h_1, \ldots, h_q)=0$ if $q\ge 2$ and $|h_1|+\cdots + |h_q|\le p$. For $|h_1|+\cdots + |h_q|>p$, $\mu(h_1, \ldots, h_q)$ can be defined recursively by noting that once $\displaystyle \mu: \bigsqcup_{1\le q'\le q-1}\mathfrak B^{q'}\to \mathbf{K}^{R}(\mathbf{x})$ has been defined, the expression on the right hand side of equation (\ref{massey}) is a cycle in homological degree $>p$; since $H_{>p}(\mathbf{K}^{R}(\mathbf{x}))=0$, we can choose $\mu(h_1, \ldots, h_q)$ such that (\ref{massey}) is satisfied.
\end{proof}

\begin{corollary}
If $I\subseteq \m^2$ is an ideal with $\mathrm{pd}_P(P/\m I)\le 3$, then $\mathfrak m I$ is Golod.  
\end{corollary}
\begin{proof}
Let $R=P/\mathfrak m I$. Due to the isomorphism $H_1(\mathbf{K}^{R}(\mathbf{x}))\cong \mathrm{Tor}_1^P(P/\m I, P/\mathfrak m)$, we can choose a basis for $H_1(\mathbf{K}^{R}(\mathbf{x}))$
consisting of elements of the form $[gdx_i]$ with $g\in I$ (corresponding to the generators $x_ig$ of $\mathfrak m I$).
 Since every cycle $z_h$ has coefficients in $\m$, it follows that $z_hz_{h'}=0$ for any $h, h'$ with $|h|+|h'|\le 3$.
 
 % Let $z_h=\sum_{|F|=i} g_Fdx_F$, then $\partial(z_h)=0$ implies $\sum_{\ell \notin G} (-1)^{\mathrm{pos}_{G\cup\{\ell\}}(\ell) }{x_{\ell}  g_{G\cup \{\ell\}}}=0$ in $P/I$, now if one of $g_{G\cup \{\ell\}}$ is unit then $\sum_{\ell \notin G} (-1)^{\mathrm{pos}_{G\cup\{\ell\}}(\ell) }{x_{\ell}  g_{G\cup \{\ell\}}}  $ is not zero in $P/\m^2$ which is a contraction.
 \end{proof}

% \begin{corollary}
%     In $P=k[x_1,x_2,x_3,x_4]$, let $I$ be a proper ideal, then $\m I$ is Golod if $\m$ is not an Associated prime of $P/I$.
% \end{corollary}

% In $P=k[x_1,x_2,x_3,x_4]$, if $\depth (P/I)=0$ and $\text{height}(I)\leq 3$, then $\m I$ is Golod?\\
% No, we have $(x_1^3,x_2^3,x_3^3,x_1x_4^2)$.\\

% In $P=k[x_1,x_2,x_3,x_4]$, if $\depth (P/I)=0$ and $\text{height}(I)\leq 2$, then $\m I$ is Golod?\\
% No, we have $(x_1^3,x_2^3,x_2x_3^2,x_1x_4^2)$.

\section{Existence of Golod Threshold with Some Upper and Lower Bounds.}
% In this section, we will introduce the notion of Golod threshold and show that it always exists; establish some general upper bound for all ideals and a lower bound for Gorenstein ideals. At the end, using these bounds, we will compute the Golod threshold for a class of Gorenstein ideals. 
Let $I$ and $J$ be two proper ideals in the standard graded polynomial ring $P=\kk[x_1,\cdots,x_n]$ with maximal homogeneous ideal $\m$ (or a regular local ring $(P,\m,\kk)$). In this section, we show that the Golod threshold of $I$ with respect to $J$ exists, and then we will focus on the case $J=\m$.
 The result below follows the steps in [Lemma 3.1, \cite{maleki2019golod}], with a slight modification.
\begin{lemma}\label{artinrees1}
    Let $M$ be a finitely generated $P$-module, and $N$ be a submodule of $M$. For any proper ideals $I$ and $J$ in $P$, there exists $d_0\in \mathbb{N}$ such that for all $d\geq d_0$ we have \[
    (J^dIM)\cap N=J^{d-d_0}(J^{d_0}IM\cap N).
    \]
\end{lemma}
\begin{proof}
    The Artin-Reese Lemma for $N\cap IM$ as a submodule of $IM$ implies that there exists $d_0$ such that for all $d\geq d_0$,\begin{equation}\label{artinrees2}
    J^dIM\cap (N\cap IM)=J^{d-d_0}(J^{d_0}IM\cap (N\cap IM))
    \end{equation}
    Observe that left hand side of \cref{artinrees2} is $J^dIM\cap N$ and the right hand side is $J^{d-d_0}(J^{d_0}IM\cap N)$.
\end{proof}

\begin{proposition}\label{torgolod}
     Let $I$ and $J$ be ideals in $P$. There exists an integer $d_0$ such that for all $d>d_0$, the natural map \begin{equation}\label{eqtorgolod}
    \Tor_i^P(P/J^dI, \kk)\xrightarrow[]{}\Tor_i^P(P/J^{d-1}I, \kk)
    \end{equation}
    induced by the surjection $P/J^dI\xrightarrow[]{}P/J ^{d-1}I$ is zero for all $i>0$.
\end{proposition}
\begin{proof}
     Set $\mathbf K^P=\mathbf{K}^P(\mathbf{x})$. Let $Z$ and $B$ be the cycles and boundaries in $\mathbf K^P$, respectively. Consider
    \[
    H_{*}(J^dI\otimes_P \mathbf K^P)\xrightarrow[]{\;\;\;\;}H_{*}(J^{d-1}I\otimes_P \mathbf K^P)
    \] 
 The map \cref{eqtorgolod} is zero if and only if \[\ker(J^dI\otimes_P \mathbf K^P_i\xrightarrow[]{}J^dI\otimes_P \mathbf K^P_{i-1})\subseteq \mathrm{im}(J^{d-1}I\otimes_P \mathbf K^P_{i+1}\xrightarrow[]{}J^{d-1}I\otimes_P \mathbf K^P_{i}),\]
 i.e., $J^dI\mathbf K^P_i\cap Z_i\subseteq J^{d-1}IZ_i$. Now, \cref{artinrees1} implies that there exists $d_0$ such that for all $d\geq d_0$, \begin{equation}\label{eq1}
 J^dI\mathbf K^P_i\cap Z_i=J^{d-d_0}(J^{d_0}I\mathbf K^P_i\cap Z_i)   
\end{equation}
for all $i\geq 0$. Since $H_i(J^{d_0}I\mathbf K^P)$ is a $\kk$-vector space we have \begin{equation}\label{eq2}
J(J^{d_0}I\mathbf K^P_i\cap Z_i)\subseteq B_{i}(J^{d_0}I\mathbf K^P)=J^{d_0}IZ_i
\end{equation}
which implies \[
 J^dI\mathbf K^P_i\cap Z_i=_{\eqref{eq1}}J^{d-d_0-1}J(J^{d_0}I\mathbf K^P_i\cap Z_i)\subseteq_{\eqref{eq2}} J^{d-d_0-1}( J^{d_0}IZ_i)=J^{d-1}IZ_i
\] 
\end{proof}

\begin{theorem}\label{existence of golod threshold}
    $\mathrm{gt}_J(I)$ exists for all proper ideals $I$ and $J$ of $P$.
\end{theorem}
\begin{proof}
    Let $L=J^{k-1}I$, where $k\geq 2$. Since $L^2\subseteq J^kI\subseteq L$, from Lemma 2.1 in \cite{maleki2019golod} and \cref{torgolod}, there exists $k_0$ such that for all $k\geq k_0$, $J^kI$ is Golod. 
\end{proof}
 We will establish upper bounds for $\mathrm{gt}_{\m}(I)$ in terms of the regularity of $I$; we recall the definiton below: 
 \begin{definition}
     Let $I\subseteq P=\kk[x_1,\ldots,x_n]$ be a proper homogeneous ideal. The regularity of $I$ denoted by $\reg(I)$ is \[
    \reg(I):= \max_{1\leq i\leq n}\{j-i\;|\; \dim_{k}[\Tor^{P}_i(P/I,k)]_j\neq 0\}.\]
 \end{definition}
The following result is stated as Theorem 6.16 (b) in \cite{mccullough2015infinite}; however, \cite{mccullough2015infinite} did not include a proof. Since we could not find a published proof in the literature, we include one here for completeness. 

\begin{proposition}\label{huneke}
    Let $P=\kk[x_1,\ldots,x_n]$ be a standard graded algebra with homogeneous maximal ideal $\m$. Let $I$ be a proper homogeneous ideal of $P$, and let $d$ be the minimal degree among the generators of $I$. If $ \reg(I)\leq 2d-3$, then $I$ is Golod. 
\end{proposition}
\begin{proof}
    Let $R=P/I$ and $\overline{\m}$ be the homogeneous maximal ideal in $R$. The assumption implies that $I\subseteq \m^2$ (since $d\ge 2$) and for $i\geq 1$, we have: \[
    \dim_{\kk}[\Tor_{i}^P(R,\kk)]_j=0 \quad\text{if}\quad j-i<d-1 \quad\text{or}\quad j-i\geq 2(d-1).
    \] 
    The conclusion follows from \cite[Lemma 9.4]{diethorn2026homological} and \cite[Proposition 9.3]{diethorn2026homological}.
\end{proof}

% We know that $\reg(\m^kI)\leq k+\reg(I)$, so if $\reg(I)\leq 2d+k-3$, then $\m^kI$ is Golod. In other words, for all $k\geq \reg(I)-2d+3$, $\m^kI$ is Golod. Therefore, an upper bound for $\mathrm{gt}_\m(I)$ is $\reg(I)-2d+3$. But we can do better; below is a general upper bound for $\mathrm{gt}_m(I)$.
% \begin{corollary}\label{upperbound}
%     If $I$ is a homogeneous ideal, we have $\mathrm{gt}_{\m}(I)\le \reg(I)-2d+3$, where $d$ is the minimal degree among the generators of $I$.
% \end{corollary}
\begin{corollary}\label{another bound}
Let $I$ be a homogeneous ideal and let $d$ and $D$ denote the smallest and respectively the largest degree of a minimal generator of $I$. Then $$\mathrm{gt}_{\m}(I) \le \mathrm{max}\left\{\dfrac{\mathrm{reg}(I)-2d+3}{2}, D-2d+2\right\}.$$
\end{corollary}
\begin{proof}
Let $k\ge \mathrm{max}\{\frac{\mathrm{reg}(I)-2d+3}{2}, D-2d+2\}$. 
Using \cite[Corollary 7.5]{eisenbud2006regularity} and $\reg(\m^k)=k-1$, we have $$\mathrm{reg}(\m^kI)\le \mathrm{max}\{\mathrm{reg}(I), k+D-1\}, 
$$
and the smallest degree of a generator of $\m^kI$ is $k+d$. The assumption implies 
$\mathrm{reg}(\m^kI)\le 2(k+d)-3$, and the conclusion follows from Proposition (\ref{huneke}).
\end{proof}

The following result will allow us to obtain another upper bound for $\mathrm{gt}_{\m}(I)$:

\begin{proposition}\label{tor zero}
Let $I$ be a homogeneous ideal and let $p:=\mathrm{pd}_P(P/I)$. If $k$ is such that 
the maps $\mathrm{Tor}_i^P(P/\m^kI, \kk)\to \mathrm{Tor}_i^P(P/I, \kk)$ induced by the canonical projection $P/\m^kI\to P/I$ are zero for all $1\le i \le \lfloor p/2\rfloor$, then $\m^kI$ is Golod. 
\end{proposition}
\begin{proof}
Consider the short exact sequence 
\begin{equation}\label{ses}
0\to \frac{I}{\m^k I} \to \frac{P}{\m^kI}\to \frac{P}{I}\to 0.
\end{equation}
 Tensoring (\ref{ses}) with ${\bf K}^P({\bf x})$ yields short exact sequences
$$
H_i({\bf K}^P({\bf x})\otimes I/\m^k I) \to H_i({\bf K}^{P/\m^kI}({\bf x}))\to H_i({\bf K}^{P/I}({\bf x})).
$$
Due to the isomorphism between Tor and Koszul homology, the rightmost map is zero for $1\le i \le \lfloor p/2 \rfloor$ by assumption. This implies that $H_i({\bf K}^{P/\m^kI}({\bf x}))$ is equal to the image of the map 
$H_i({\bf K}^P({\bf x})\otimes I/\m^kI)\to H_i({\bf K}^{P/\m^kI}({\bf x}))$. 

    % ({\em note: commenting out the thing about coefficients in $\m^{k-1}I$})
    % We claim that the assumption that the ideal is equigenerated implies that all such cycles must in fact have coefficients in $\m^{k-1}I$\dinesh{This also means for $i\leq \lfloor p/2 \rfloor $ we have monomial basis for $H_i({\bf K}^P\otimes P/\m^kI) $}. Let $d$ denote the degree of the generators of $I$. An element in $H_i({\bf K}^P\otimes P/m^kI)$ must have degree at least $d+k+i-1$. If $[z]=\sum_{|F|=i} g_F dx_F$ is homogeneous of degree $\ge d+k+i-1$, then $\mathrm{deg}(g_F)\ge d+k-1$ and, given $g_F\in I$, this implies $g_F\in \m^{k-1}I$.

It follows that for all $1\le i\le \lfloor p/2 \rfloor$, every element of $H_i({\bf K}^{P/\m^kI}({\bf x}))$ can be represented by a cycle that has coefficients in $I/\m^kI$.
Also note that for all $i\ge 1$, every nonzero element in $H_i({\bf K}^{P/\m^kI}({\bf x}))$ must have degree at least $k+i$, and therefore it is represented by a cycle with coefficients in $\m^k/\m^kI$. 
    
    Consider a basis $\{[z_{\lambda}]\}$ for $H_{\ge 1}({\bf K}^P({\bf x})\otimes P/\m^kI)$ with $[z_{\lambda}]\in H_{i_{\lambda}}$. 
    Let $\lambda_1, \lambda_2$ be such that $i_{\lambda_1}+i_{\lambda_2}\le p$. 
 We may assume without loss of generality that $i_{\lambda_1}\le \lfloor \frac{p}{2}\rfloor.$
 Since $z_{\lambda_1}$ has coefficients in $I/\m^kI$ and $z_{\lambda_2}$ has coefficients in $\m^k/\m^kI$, it follows that $z_{\lambda_1}z_{\lambda_2}=0$ and the conclusion follows from
\cref{trivialmassey1}.
\end{proof}
\begin{corollary}\label{upper bound}
   Let $I$ be a homogeneous ideal with $p:=\mathrm{pd}_P(P/I)$, and $d$ be the smallest degree of a generator of $I$.
For each $i\le p$, let $$T(i):=\mathrm{max}\{j: [\Tor_{i}^P(P/I,\kk)]_j\ne 0\}.$$
   Then
$$
\mathrm{gt}_{\m}(I)\le \mathrm{max}\{T(i)-i-d+2: i \le \lfloor p/2 \rfloor \}
$$ 
\end{corollary}
\begin{proof}
Let $k\ge \mathrm{max}\{T(i)-i-d+2: i \le \lfloor p/2 \rfloor \}$.
  We have $[\mathrm{Tor}_i(P/I, \kk)]_j=0$ for $j>T(i)$, and 
  $[\mathrm{Tor}_i(P/\m^kI, \kk)]_j=0$ for $j<d+k+i-1.$
  The assumption implies that for all $i\le \lfloor \frac{p}{2}\rfloor$, we have $d+k+i-1>T(i)$, and therefore the map
  $$
[\mathrm{Tor}_i(P/I, \kk)]_j\to [\mathrm{Tor}_i(P/\m^kI, \kk)]_j
  $$
  is zero in each degree $j$. \cref{tor zero} shows that $\m^kI$ is Golod and the conclusion follows.
\end{proof}

\begin{example}
Let $I=(x_1^d,\ldots,x_c^d)$ be a monomial complete intersection in $P=\kk[x_{1},\cdots,x_{n}]$. When $c$ is even, both \cref{another bound} and \cref{upper bound} give the same upper bound for $\mathrm{gt}_m(I)$, namely $(\frac{c}{2}-1)(d-1)+1$; when $c$ is odd, \cref{upper bound} gives a better upper bound, namely $(\lfloor \frac{c}{2}\rfloor -1)(d-1)+1$. It turns out that the value of $\mathrm{gt}_{\m}(I)$ is in fact equal to this upper bound. We will calculate the value of $\mathrm{gt}_{\m}(I)$ for a general monomial complete intersection $I=(x_1^{a_1}, \ldots, x_n^{a_n})$ in \cref{preciseCIbound}.

% We have $\mathrm{pd}_P(P/I)=c$,  $T(i)=id$ for all $1\le i \le c$ and $\reg(I)=c(d-1) $. Suppose $d\geq 2$, then for $i\leq c$, we have \[
% \mathrm{max}\{T(i)-i-d+2: i \le \lfloor c/2 \rfloor \}=\left({\left\lfloor \frac{c}{2}\right\rfloor}-1\right)(d-1)+1
% \]
% and 
% \[
% \mathrm{max}\left\{\dfrac{\mathrm{reg}(I)-2d+3}{2}, D-2d+2\right\}=\left(\frac{c}{2}-1\right)(d-1)+1
% \]
% Therefore, when $c$ is even both \cref{another bound} and \cref{upper bound} give the same upper bound for $\mathrm{gt}_m(I)$, and when $c$ is odd \cref{upper bound} gives a better upper bound.

\end{example}

If the ideal $I$ is Gorenstein, the following provides a lower bound for $\mathrm{gt}_{\m}(I)$:

\begin{theorem}\label{lower bound for gorenstein}
Let $I\subseteq P$ be a homogeneous Gorenstein ideal of codimension $c$ and regularity $\mathrm{reg}(I)$. Let $b_{ij}$ denote the dimension of $[\mathrm{Tor}_i(P/I, \kk)]_j$. For $0\le i\le c$, let $t(i)$ and $T(i)$ denote the smallest and respectively largest values of $j$ such 
$b_{ij}\ne 0$ and let $D$ be the largest degree of a generator of $I$.

Let $(i_1, j_1)$ be such that $1\le i_1\le c-1$ and $b_{i_1 j_1}\ne 0$, and let $i_2=c-i_1$, $j_2=c+\mathrm{reg}(I)-j_1$.

Then 
$$
\mathrm{gt}_{\m}(I)> \mathrm{min}\{j_1-i_1-D+1, j_2-i_2-D+1\}
$$\end{theorem}
Before proving this theorem, we will need a Lemma:
\begin{lemma}\label{mod k}
Let $R=P/I$ and $R_k=P/\m^kI$ where $I$ is a homogeneous ideal of the polynomial ring $P$, and let $D$ denote the largest degree of a generator of $I$. Let $\pi_k: {\bf K}^{R_k}\to {\bf K}^R$ denote the map of complexes induced by  canonical projection $\pi: R\to R_k$.
For every $[z]\in H_i({\bf K}^{R})$ with
$\mathrm{deg}([z])\ge k+D+i-1$, there exists $[\tilde{z}]\in H_i({\bf K}^{R_k})$ such that $[z]=\pi_k(\tilde{[z]})$. 
\end{lemma}
\begin{proof}
Write $$[z]=[\sum_{|F|=i} \overline{g}_Fdx_F]$$
where the summation is over subsets $F\subseteq [n]$ and $g_F$ are elements of $P$ with $\overline{g}_F$ being their image in the quotient ring $R$. Let $\tilde{g}_F\in R_k$ be the image of $ g_F$ in the quotient ring $R_k$.
The assumption that $[z]\in H_i({\bf K}^R)$ means that $\sum_{|F|=i} \overline{g}_Fdx_F$ is a cycle in the Koszul complex, i.e. for each $G\subseteq [n]$ with $|G|=i-1$ we have 
$$
\sum_{\ell \notin G} (-1)^{\mathrm{pos}_{G\cup\{\ell\}}(\ell) }\overline{x_{\ell } g_{G\cup \{\ell\}}}=0,$$ 
and therefore 
$$
\sum_{\ell \notin G} (-1)^{\mathrm{pos}_{G\cup\{\ell\}}(\ell) }x_{\ell } g_{G\cup \{\ell\}}\in I.$$
The degree assumption implies that 
$$
\sum_{\ell \notin G} (-1)^{\mathrm{pos}_{G\cup\{\ell\}}(\ell)} x_{\ell } g_{G\cup \{\ell\}}\in \m^k I,$$
and therefore 
$\tilde{z}:=\sum_{|F|=i} \tilde{g}_F dx_F$ is a cycle in ${\bf K}^{R_k}$.
\end{proof}

We are now ready to prove \cref{lower bound for gorenstein}:
\begin{proof}
Let
\begin{equation}\label{lower bound}
k\le \mathrm{min}\{j_1-i_1-D+1, j_2-i_2-D+1\}
    \end{equation}
% $0\le i \le c$ 
% such that 
% \begin{equation}\label{lower bound2}
% k\le \mathrm{min}\{t(i)-i-D+1, \mathrm{reg}(I)-T(i)+i-D+1\}
% \end{equation}
We will prove that $\m^k I$ is not Golod.

Due to the isomorphism $\mathrm{Tor}_i^P(P/I, \kk)\cong H_i({\bf K}^R)$, $b_{i_1 j_1}\ne 0$ implies that there is a nonzero $[z_1]\in H_{i_1}({\bf K}^{R})$ with $\mathrm{deg}([z_1])=j_1$, and Poincar\'{e} duality implies that there exists a complementary $[z_2]\in H_{i_2}({\bf K}^R)$ with $\mathrm{deg}([z_2])=j_2$ such that $[z_1][z_2]\ne 0$ in $H_c({\bf K}^R).$
\cref{mod k}
yields
$[\tilde{z_1}]\in H_{i_1}({\bf K}^{R_k})$ with $\pi_k([\tilde{z_1}])=[z_1]$ and 
$[\tilde{z_2}]\in H_{i_2}({\bf K}^{R_k})$ with $\pi_k([\tilde{z_2}])=[z_2]$. Since $[z_1][z_2]\ne 0$, it follows that $\tilde{[z_1]}\tilde{[z_2]}\ne 0$, and therefore $\m^k I$ is not Golod.
\end{proof}

% {\em commenting out the way proof was originally written }

% Since $j\ge t(i)$, the assumption (\ref {lower bound}) implies $\mathrm{deg}([z])\ge k+D+i-1$, and Observation (\ref{mod k}) implies that there exists $[\tilde{z_1}]\in H_i(\tilde{\bf K}_k)$ with $\pi_k([\tilde{z_1}])=[z_1]$.   

% Poincar\'{e} duality implies that there exists $[z_2]\in H_{c-i}({\bf K})$ such that $[z_1][z_2]\ne 0$, and $\mathrm{deg}([z_2])=c+\mathrm{reg}(I)-j$. Since $j\le T(i)$, the assumption (\ref{lower bound}) implies that $\mathrm{deg}([z_2])\ge k+D+(c-i)-1$, and Observation (\ref{mod k}) yields a $\tilde{[z_2]}\in H_{c-i}(\tilde{\bf K}_k)$ such that $[z]_2=\pi_k(\tilde{[z_2]})$. Since $[z_1][z_2]\ne 0$, it follows that $\tilde{[z_1]}\tilde{[z_2]}\ne 0$, and therefore $\m^k I$ is not Golod.
\begin{example}\label{cyclicGorensteinexample}
For $n\geq 5$ and $a\geq 1$, consider the ideal 
$$
I_{a,n}=(x_i^ax_j^a:\; \ i, j \in [n] \ \text{and} \ |i-j|\neq 0,1,n-1)\subset P=\kk [x_1, \ldots, x_n]
$$
Then 
$$
\mathrm{gt}_{\m}(I_{a, n})=\left\lfloor \frac{n-4}{2}\right\rfloor  (a-1)+1
$$
\end{example}
\begin{proof} We will describe the pairs $(i, j)$ for which $b_{ij}\ne 0$, where $b_{ij}$ are graded betti numbers of $P/I_{a, n}$, and then apply \cref{upper bound} and \cref{lower bound for gorenstein}.
We start by observing that $P/I_{1, n}$ is the Stanley-Reisner ring of the cycle $C_n$ with vertices $\{1, \ldots, n\}$.
For a squarefree monomial $\displaystyle M=\prod_{i\in W} x_i$, Hochster's formula \cite[Theorem 5.1]{hochster1977cohen} for the multigraded Betti numbers gives
$$
b_{i, M}(P/I_{1, n})=\mathrm{dim}_{\kk} (\tilde{H}_{|W|-i-1}(\Delta_W, \kk))
$$
where $\Delta_W$ is the simplicial complex induced by $C_n$ on $W$.

For $|W|<n$, $\Delta_W$ is a disjoint union of paths and therefore 
$\tilde{H}_r(\Delta_W)=0$ for all $r>0$, and $\tilde{H}_0(\Delta_W)=c(W)-1$, where $c(W)$ denotes the number of connected components of $\Delta_W$. In particular, 
$
b_{i, M}(P/I_{1, n})\ne 0$ if and only if $i=\mathrm{deg}(M)-1$ and $\Delta_W$ is disconnected, i.e. $M$ is not a product of $i+1$ consecutive variables (where $x_n, x_1$ are considered consecutive).

For $W=[n]$ (i.e. $M=\prod_{i=1}^n x_i$), we have $\Delta_W=C_n$, and $\tilde{H}_r(C_n)\ne 0\Leftrightarrow r=1$; more precisely, $\tilde{H}_1(C_n)\cong \kk$. It follows that 
$b_{i, M}\ne 0\Leftrightarrow i=n-2$. 
We see that $\mathrm{pd}_P(P/I_{1, n})=n-2$ and
$b_{n-2, M}\ne 0\Leftrightarrow M=\prod_{i=1}^n x_i$, in which case $b_{n-2, M}=1$.

In view of Theorem 3.3 in \cite{gasharov1999lcm}, the multigraded Betti numbers in the resolution of $I_{a, n}$ are obtained from those in the resolution $I_{1, n}$ by replacing each monomial $M=\prod_{i\in W} x_i$ by $M_a:=\prod_{i\in W} x_i^a$. In terms of the standard grading, we have $b_{i, j}\ne 0 \Leftrightarrow j=(i+1)a$ for $i\le n-3$, and $b_{n-2, j}\ne 0\Leftrightarrow j=na$.

We apply \cref{upper bound} to obtain an upper bound for $\mathrm{gt}_{\m}(I_{a, n})$:
\begin{equation}\label{gorenstein example}
    \begin{split}
\mathrm{gt}_{\m}(I_{a, n})& \le \mathrm{max}\{ (i+1)a-i-2a+2  \, : \,   i\le \left\lfloor \frac{n-2}{2}\right\rfloor \}\\
&=\left\lfloor \frac{n-4}{2}\right\rfloor (a-1)+1\\
\end{split}
\end{equation}
Consider the ring homomorphism $\psi:P\xrightarrow[]{}P$, where $\psi(x_i)=x_i^a$. Here $\psi$ is faithfully flat, as $P$ is a free $P$-module via $\psi$. Changing the base from $P$ to $P/I_{1,n}$, we have the flat map $ \psi\otimes P/I_{1,n} :P/I_{1,n}\xrightarrow[]{}P/I_{a,n}.$ Now, since $P/I_{1,n}$ and $P/(x_1^a,\ldots,x_n^{a})$ are Gorenstein, after localizing at $\m$ and using \cite[Corollary 3.3.15]{bruns1998cohen}, we have $P/I_{a,n}$ is Gorenstein.

 We can now use \cref{lower bound for gorenstein} to see that have equality in equation \cref{gorenstein example}.
Indeed, let $(i_1, j_1)$, $(i_2, j_2)$ be as in the assumption of \cref{lower bound for gorenstein}, so that $i_1+i_2=n-2$, $j_1=(i_1+1)a$, $j_2=(i_2+1)a$. Without loss of generality, assume $i_1\le i_2$. We have
$$
\mathrm{gt}(I_{a, n})>\mathrm{min}\{(i_1-1)(a-1), (i_2-1)(a-1)\}=(i_1-1)(a-1)
$$
Since this holds for any choice of $i_1, i_2$ subject to $i_1+i_2=n-2$, we can take $i_1=\lfloor \frac{n-2}{2}\rfloor$, $i_2=\lceil \frac{n-2}{2}\rceil$ and the conclusion follows.
\end{proof}

The proof of \cref{cyclicGorensteinexample} shows that for the Gorenstein ideal $I_{a,n}$, if $m^kI_{a,n}$ is Golod then $m^{k'}I_{a,n}$ is also Golod for all $k'>k$. The following example shows that this behavior does not always hold.

\begin{example}
    Let $I=(x_1^3, x_2^3, x_3^3,x_4^3, x_1x_2^2x_3^2x_4^2)\subseteq \kk[x_1, \ldots, x_4]$. For this ideal, $\m I$ is Golod (one can check that the criterion \cite[Theorem 1.5]{zoromski2024monomial} is satisfied for $\m I$) but $\m^2I$ is not.  Indeed,    $[x_1^2x_2^2dx_{12}]$ and $[x_3^2x_4^2dx_{34}]$ are elements in $H_{\ge 1}(K^{P/\m^2I}(\mathbf x))$
    and their product is not zero. 
\end{example}

The above example also resolves Question 5.5 \cite{de2016products} in the negative, as $\m I $ and $\m$ are Golod, but $\m\cdot \m I$ is not Golod. Now, we note that even though there are examples where the Golod property of $\m^k I$ oscillates, there are an abundance of ideals for which the Golod threshold is sharp, i.e., if $\m^kI$ is Golod then $\m^{k'}I$ is Golod for all $k'\geq k$. To prove this, we will show that we can pull back the trivial Massey operation under a suitable degree condition. For a similar construction under different hypotheses, one can look at \cite{vandebogert2022detecting}.
 
 \begin{theorem}\label{equigenratedgolodThreshold}
     Let $I \subseteq \m^2$ be a homogeneous ideal in a standard graded polynomial ring $P$. Let $d$ and $D$ denote the smallest and largest degrees, respectively, of a minimal generator of $I$. If $D < 2d-1$ and $I$ is Golod, then $\m I$ is also Golod. 
 \end{theorem}
 
 \begin{proof}
 Consider the canonical short exact sequence
     \[
     0\xrightarrow[]{}\mathbf{K}^{I/\m I}(\textbf{x})\xrightarrow[]{} \mathbf{K}^{P/\m I}(\textbf{x})\xrightarrow[]{\varphi}\mathbf{K}^{P/I}(\textbf{x})\xrightarrow[]{}0,
     \]
     and the exact sequence of homology \[
      H(\mathbf{K}^{I/\m I}(\textbf{x}))\xrightarrow[]{} H(\mathbf{K}^{P/\m I}(\textbf{x}))\xrightarrow[]{\varphi_*}H(\mathbf{K}^{P/I}(\textbf{x}))
     \]
     Observe that 
    $ [\mathbf{K}_i^{I/\m I}(\textbf{x})]_j=0$ for $j>i+D$. Therefore, if 
    $z\in [\mathbf{K}_i^{P/I}(\textbf{x})]_j$
    with $j>i+D$, there is a unique homogeneous preimage $\varphi^{-1}(z)\in [\mathbf{K}_i^{P/\m I}(\textbf{x})]_j.$
    Let $\mathfrak B$ be a homogeneous basis for $ H_{\geq 1}(\mathbf{K}^{P/ I}(\textbf{x}))$ with trivial Massey operation $\mu$ on $\mathfrak B$ and let $\mathfrak B'$ be a homogeneous basis for $ H_{\geq 1}(\mathbf{K}^{P/\m I}(\textbf{x}))$. We will use $\mu$ to define a trivial Massey operation $\mu'$ on $\mathfrak B'$. First, extend $\mu$ multilinearly to $\bigsqcup_{q\geq 1}H_{\geq 1}(\mathbf{K}^{P/I}(\textbf{x}))^{q}$; note that equation (\ref{massey}) continues to hold for all $h_1, \ldots, h_q\in H_{\geq 1}(\mathbf{K}^{P/I}(\textbf{x}))$; see \cite[Lemma 11.2]{avramov2013contravariant}.
    For $h\in \mathfrak B'$, if $\varphi_*(h)\neq 0$, we can pick a representative $z_h$
    such that $\mu(\varphi_*(h))=\varphi (z_h)$, define $\mu'(h)=z_h $ and if $\varphi_*(h)=0$, we can pick a representative $z_h$ such that $z_h\in I\mathbf{K}^{P/\m I}(\textbf{x})$, define $\mu'(h)=z_h $. For $(h_1,\ldots,h_{r})\in \bigsqcup _{q\geq 2} \mathfrak B'^q$, if $\varphi_*(h_i)=0$ for some $i$, define $\mu'(h_1,\ldots,h_{r})=0$. If  $\varphi_*(h_i)\neq 0$, for all $i$, let $i_{h_t}$ and $j_{h_t}$ denote the homological and respectively the internal degree of $h_t$ for each $t$. We have  
    \[\mu(\varphi_{*} (h_1),\ldots,\varphi_{*}(h_r))\in [\mathbf K_i^{P/I}({\bf x})]_j, 
    \text{with} \]
    \[i=\sum_{t=1}^r i_{h_t}+r-1 \ \text{and}\ j=\sum_{t=1}^r j_{h_t}.\] Note that we have $j_{h_t}\ge i_{h_t}+d$
    for all $t$, and therefore $j\ge i+rd-r+1\ge i+D+1.$
     Define $ \mu'(h_1,\ldots,h_r): = \varphi^{-1} (\mu(\varphi_{*} (h_1),\ldots,\varphi_{*}(h_r))\in [\mathbf{K}_i^{P/\m I}(\textbf{x})]_j$ (as observed earlier, there is a unique such element).

    It is clear that $\mu'$ satisfies conditions $(i)$ and $(iii)$ in the definition of a trivial Massey operation. It remains to check condition $(ii)$. Let $(h_1,h_2)\in \mathfrak B'^2$. If $\varphi_*(h_1) $ or $\varphi_*(h_2)$ is zero then one of $z_{h_1}$ or $z_{h_2}$ has coefficients in $I$, and the other one has coefficients in $\m$; it follows that $\overline{z_{h_1}}\cdot z_{h_2}=0=\partial (\mu'(h_1, h_2))$. If both $\varphi_*(h_1) $ and $\varphi_*(h_2)$ are not zero, then \begin{equation*}
    \begin{aligned}
    \partial(\mu'(h_1,h_2))&=\partial(\varphi^{-1}(\mu(\varphi_*(h_1),\varphi_*(h_2))))\\
    &=\varphi^{-1}(\partial\mu(\varphi_*(h_1),\varphi_*(h_2)))\\
    &=\varphi^{-1}(\overline{\mu(\varphi_*(h_1))}\cdot \mu(\varphi_*(h_2)))\\
    &=\varphi^{-1}(\overline{\varphi(z_{h_1})}\cdot\varphi(z_{h_2}))\\
    &=\overline{z_{h_1}}\cdot z_{h_2}\\
    &=\overline{\mu'(h_1)}\mu'(h_2).  
    \end{aligned}
    \end{equation*}
    Now, let $q>2$. If $ \varphi_*(h_i)=0$ for some $i$, then  $ \sum_{j=1}^{q-1}\overline{\mu'(h_1,\cdots,h_{j})}\mu'(h_{j+1},\cdots,h_q)=0=\partial(\mu'(h_1, \cdots, h_r))$. If $ \varphi_*(h_i)$ is not zero for all $i$, then\begin{equation*}
    \begin{aligned}
    \partial(\mu'(h_1,\ldots,h_q))&=\partial(\varphi^{-1}(\mu(\varphi_*(h_1),\ldots,\varphi_*(h_q))))\\
    &=\varphi^{-1}(\partial((\mu(\varphi_*(h_1),\ldots,\varphi_*(h_q)))))\\
    &=\varphi^{-1}(  \sum_{j=1}^{q-1}\overline{\mu(\varphi_*(h_1),\cdots,\varphi_*(h_{j}))}\mu(\varphi_*(h_{j+1}),\cdots,\varphi_*(h_r)))\\
    &=  \sum_{j=1}^{q-1}\overline{\varphi^{-1}\mu(\varphi_*(h_1),\cdots,\varphi_*(h_{j}))}\varphi^{-1}\mu(\varphi_*(h_{j+1}),\cdots,\varphi_*(h_r))\\
    &=  \sum_{j=1}^{q-1}\overline{\mu'(h_1,\cdots,h_{j})}\mu'(h_{j+1},\cdots,h_r)
    \end{aligned}
    \end{equation*}
\end{proof}

\begin{corollary}
    Let $(P,\m,\kk)$ be a standard graded polynomial ring and $I\subseteq \m^2$ a proper homogeneous ideal. Let $d$ and $D$ be the lowest and respectively the highest degree of a minimal generator of $I$. If $D < 2d-1$, then $ \m^{k}I$ is Golod if and only if $k\geq\mathrm{gt}_{\m}(I)_{}$. 
\end{corollary}

\section{Combinatorial Criterion for Monomial Koszul Basis}
% In this section we give a combinatorial criterion for the existence of a monomial basis of Koszul homology of a monomial quotient ring. Using this criterion, we will give a sufficient and necessary condition for monomial ideals in 4 variables to have a monomial basis for their Koszul Homology. At the end end we will show that the ideal of the form $ \m^k (x_1^{a_1},\cdots,x_n^{a_n})$ admits a monomial basis for its Koszul homology, and we will use this result in next section to establish the exact formula for the Golod threshold for a monomial complete intersection.

Let $I$ be a monomial ideal in $P=\kk[x_1,\cdots,x_n]$ and $R=P/I$. Fix $\alpha\in \mathbb{N}^n$ and let $V=\mathrm{Supp}(x^{\alpha})$. The Lower Koszul complex for $I$ in degree $\alpha$ is  the simplicial complex \[\mathbf{K}_{\alpha}(I)=\{G\subseteq V:\; x^{\alpha-1_{V\setminus G}}\notin I\}.\]  $\mathbf{K}_{\alpha}(I)$ consists of the complements of the upward-closed set $\mathfrak T_{\alpha}= \{F\subseteq V: x^{\alpha-1_F}\notin I \}$, and therefore it is a simplicial complex. Note that the  elements $F\in\mathfrak T_{\alpha}$ are in one-to-one correspondence with the basis elements $x^{\alpha-1_{F}}dx_F$ of the $\kk$ vector space $[\mathbf{K}^{R}(\mathbf{x})]_{\alpha}$.
\begin{definition}
A {\em monomial cycle} in $\mathbf{K}^{R}(\mathbf{x})$ is a cycle of the form $[x^{\alpha-1_F}dx_F]$
for some $F\subseteq [n]$. 
We say that the Koszul homology of $R$ admits a {\em monomial Koszul basis} if there exists a basis of $H_{\ge 1}(\mathbf{K}^{R}(\mathbf{x}))$ consisting of monomial cycles. 
    
\end{definition}

We will now give a sufficient condition for the Koszul homology of $R$ to admit a monomial Koszul basis, which will be in terms of shellability of the lower Koszul complexes. For the definition of shellability, we refer the reader to \cite{bjorner1996shellable}.

% The following Lemma shows that the monomial cycles in $[\mathbf{K}^{R}(\mathbf{x})]_{\alpha}$ come  from the facets of $\mathbf{K}_{\alpha}(I)$.  

% \begin{lemma}\label{Facetcycle}
% Let $F\subseteq V$ and $G=V\setminus F$, suppose
% \(x^{\alpha-1_F}\notin I\).  Then $x^{\alpha-1_F}dx_F$ is a cycle in $[\mathbf{K}^R(\textbf{x})]_{\alpha}$ if and only if $G$ is a facet in $\mathbf K_\alpha(I)$.
% \end{lemma}

% \begin{proof}
% $\partial(x^{\alpha-1_F}dx_F)=0\iff
% \forall i\in F,\;\;x^{\alpha-1_{F\setminus\{i\}}}\in I
% \iff\forall i\notin G,\;\;
% G\cup\{i\}\notin \mathbf K_\alpha(I)\iff
% G \text{ is a facet in }\mathbf K_\alpha(I)$.
% \end{proof} 

\begin{theorem}
    \label{shellmonbasis}
     If $K_\alpha(I)$ is shellable for all multidegrees $\alpha$,  then $H_{\geq 1}(\mathrm{ K}^R(\textbf{x}))$ admits a monomial basis.
\end{theorem}
\begin{proof}
     For each multidegree $\alpha\in \mathbb{N}^n$, consider the chain isomorphism,
     \begin{equation}\label{chain map}
\Psi:[\mathbf K^R(\textbf{x})]_{\alpha}\xrightarrow[]{\;\;\;\;\;\;}\widetilde{\mathbf{C}}^{\cdot}(\mathbf{K}_{\alpha}(I);\kk)
\end{equation}
\[
x^{\alpha-1_F}dx_F\xrightarrow[]{\;\;\;\;\;}(-1)^{\sigma(F)+|F|}(V\setminus F)^*
\]
where $\sigma(F)=\sum_{i\in F}\mathrm{pos}_{V}(i)$. Since $\mathbf{K}_{\alpha}(I)$ is shellable, a basis for $H^{\geq 1}(\widetilde{\mathbf{C}}^{\cdot}(\mathbf{K}_{\alpha}(I);\kk))$ exists, whose class representatives can be chosen to be dual to facets in $\mathbf{K}_{\alpha}(I)$ by \cite[Theorem 4.3]{bjorner1996shellable}. From the isomorphism $\Psi$, we get a corresponding basis for $H_{\geq 1}([\mathbf K^R(\textbf{x})]_{\alpha})$, which is monomial by construction. The union of all the multigraded bases gives a monomial basis for $H_{\ge 1}(\mathbf K^R(\textbf{x})).$
\end{proof}
% As a result of [\cite{gasharov1999lcm}, Lemma 2.2], to check the monomial basis for $H_{\geq 1}(\mathbf{K}^R(\textbf{x}))$, we just need to check that the lower Koszul complex for $I$ at multidegree $\alpha$ is shellable for only finitely many $\alpha$. 

One immediate consequence of \cref{shellmonbasis} is that a monomial ideal in at most 3 variables admits a monomial basis for its Koszul Homology; note that this was already known from \cite{dao2022monomial}. We note that the converse of \cref{shellmonbasis} holds in 4 variables.

\begin{proposition}\label{4varclassification}
    Let $I$ be a monomial ideal in $\kk[x_1,x_2,x_3,x_4]$. The Koszul homology for $I$ admits a monomial basis if and only if all the lower Koszul complexes for I at every multidegree are shellable. 
 %    \begin{equation}\label{4shell}
 %    x^{\alpha-1_{\{1,2\}}},\;x^{\alpha-1_{\{3,4\}}}\notin I\quad \text{and}\quad
 % x^{\alpha-1_{\{1,3\}}},\;x^{\alpha-1_{\{1,4\}}},\;x^{\alpha-1_{\{2,3\}}},\;x^{\alpha-1_{\{2,4\}}}\in I        
 %     \end{equation}
 \end{proposition}
\begin{proof}
     Suppose there exists $\alpha=(\alpha_1,\cdots,\alpha_4)\in \mathbb{N}^4$, such that $ \mathbf{K}_{\alpha}(I)$ is not shellable. Up to labeling, there is a unique non-shellable simplicial complex on 4 vertices, which is two disjoint one-dimensional facets. Therefore we can assume $ \mathbf{K}_{\alpha}(I)=\{ \{1,2\},\{3,4\} \} $.  Consider the cycle $z=x^{\alpha-1_{\{2,3,4\}}}dx_{234}-x^{\alpha-1_{\{1,3,4\}}}dx_{134}$. The boundary of $[\mathbf K^R(\textbf{x})]_{\alpha}$ at homological degree 3 is generated by the differential of $x^{\alpha-1_{\{1,2,3,4\}}}dx_{1234}$; therefore, the class $[z]$ is non-zero, and since there are no monomial cycles in $[\mathbf K^R(\textbf{x})]_{\alpha}$ at homological degree 3, the class $[z]$ cannot be represented by a monomial basis. Therefore, the Koszul homology does not admit a monomial basis. 
The other implication is known from \cref{shellmonbasis}. 
 \end{proof}

% \begin{example}
%     Fix an integer $k\geq 1$. Let $I_{k}=(x_1^{k},x_2^{k})(x_3,x_4)$. For $\alpha=(k,k,1,1)$, the lower Koszul complex for $I_{k}$ is not shellable; therefore, the Koszul homology for $I_k$ does not have a monomial basis. As one can check, \[H_{3}([\mathbf{K}^R(\textbf{x})]_{(k,k,1,1)})\cong \kk[x_1^kx_2^{k-1}dx_{234}-x_1^{k-1}x_2^{k}dx_{134}].
%     \]
% \end{example}

Throughout the rest of this section, let $P=k[x_1,\ldots,x_n]$ be a standard multigraded polynomial ring with homogeneous maximal ideal $\mathfrak m=(x_1,\ldots,x_n)$, and let
\[
I=(x_1^{a_1},\ldots,x_n^{a_n}), \quad a_i\geq 1 \text{ for all } 1\leq i\leq n,
\]
be a monomial complete intersection. Set $R=P/I$. We will use \cref{shellmonbasis} to show that the Koszul homology of $\m^k I$ admits a monomial basis for all $k\geq 0$.

Fix a multidegree $\alpha\in \mathbb{N}^n$. For each $i\in [n]$, let $q_i:=|\alpha| -a_i-k$ and $q'_i:=|\mathrm{Supp}(x^{\alpha}) | - q_i$. Recall that the lower Koszul complex $\mathbf{K}_{\alpha}(\m^kI)$ is the complement of the upward-closed set \[
\mathfrak T_{\alpha}=\{F\subseteq \mathrm{Supp}(x^{\alpha}): x^{\alpha-1_F}\notin \m^kI \}
\]

\begin{obs}\label{obs1} Let $F\subseteq \mathrm{Supp}(x^{\alpha})$. We have $F\in \mathfrak T_{\alpha}$ if and only if for all $i\in [n]$ with $\alpha_i\ge a_i$, we have ($\alpha_i=a_i$ and $i\in F$) or $|F|>q_i$.

\end{obs} 

\begin{obs}\label{obs2}
For $G\subseteq \mathrm{Supp}(x^{\alpha})$, we have $G\in \mathbf K_{\alpha}(I)$ if and only if for all $i$ with $\alpha_i\ge a_i$, we have ($\alpha_i=a_i$ and $i\notin G$) or $|G|<q_i'$. Equivalently, $G\in \mathbf K_{\alpha}(I)$ if and only if the following hold:

\begin{enumerate}
\item $|G|<\mathrm{min}\{q_i': i\in [n]\ \mathrm{with}\ \alpha_i>a_i\}$, and 

\item $|G|<\mathrm{min}\{q_i': i\in G \ \mathrm{with}\ \alpha_i=a_i\}$

\end{enumerate}

\end{obs}
\begin{proof}
Apply Observation \ref{obs1} to $F:=\mathrm{Supp}(x^{\alpha})\, \backslash\, G$.
\end{proof}

\begin{proposition}\label{shiftedsimplex}
There is a suitable ordering of the variables $x_1, \ldots, x_n$ such that the simplicial complex $\mathbf K_{\alpha}(I)$ is shifted, i.e. if $G\in \mathbf K_{\alpha}(I)$, $j\in G$ and $j'\notin G$ such that $j'<j$, then $G':=(G\, \backslash \{j\}) \cup \{j'\} \in \mathbf K_{\alpha}$.
\end{proposition}
\begin{proof}
Reorder the variables $x_1, \ldots, x_n$ in such a way that the indices $i$ with $\alpha_i\ne a_i$ (if any) are in the beginning of the list, i.e. 
there exists $s\in [n]$ such that $\alpha_i=a_i \Leftrightarrow i\ge s+1$, and $q'_{s+1}\ge q'_{s+2}\ge \ldots \ge q'_n$.

Let $G\in \mathbf K_{\alpha}(I)$, $j\in G$ and$j'\notin G$ with $j'<j$, and let $G':=(G\, \backslash \{j\}) \cup \{j'\}$. We need to prove that $G'$ satisfies conditions (1) and (2) from Observation (\ref{obs2}). Note that $|G'|=|G|$, so condition (1) for $G'$ follows immediately from condition (1) for $G$. The ordering of the variables together with the assumption $j'<j$ implies that 
$$
\mathrm{min}\{q_i': i \in G \ \mathrm{with} \ \alpha_i=a_i\} \le \mathrm{min}\{q_i': i \in G'\  \mathrm{with} \ \alpha_i=a_i\},
$$
and condition (2) for $G'$ also follows from the corresponding condition for $G$.
\end{proof}
\begin{corollary}\label{monomialbasisforCI}
    The Koszul homology of  $\m^k I$ admits a monomial basis for all $k\geq 0$.
\end{corollary}
\begin{proof}
    From \cref{shiftedsimplex} we know that the lower Koszul complex for the ideal $\m^k I$ is shifted at all multidegrees. Since all shifted complexes are shellable \cite[Theorem 11.3]{bjorner1997shellable}, \cref{shellmonbasis} implies that the Koszul homology of $\m^k I$ admits a monomial basis.
\end{proof}

\begin{remark}
    Using a similar argument, one can also show that the conclusion of \cref{monomialbasisforCI} also holds for ideals of the form $\sum_{i=1}^nx_i^{a_i}\m^{c_i}$. 
\end{remark}

\section{The Golod Threshold for monomial complete intersection ideals}
% In this section, we will compute the Golod threshold for a monomial complete intersection. 
Throughout this section, let $P=k[x_1,\ldots,x_n]$ be a standard graded polynomial ring with homogeneous maximal ideal $\mathfrak m=(x_1,\ldots,x_n)$, and let
\[
I=(x_1^{a_1},\ldots,x_n^{a_n}), \quad a_i\geq 1 \text{ for all } 1\leq i\leq n,
\]
be a monomial complete intersection. Set $R=P/I$. 

 For $F\subseteq [n]$ with $|F|\geq 2$, let \[
a_F:=\sum_{i\in K} a_i-|F|=\sum_{i\in F}(a_i-1)
\]
and \[
m_F:=\mathrm{min}\{a_{F'}\, | \, F'\subseteq F, |F'|=|F|-1\}
\]
For a partition $K\sqcup L= [n]$ with both $K$ and $L$ consisting of at least two elements, let $m_{K, L}:=\mathrm{min}\{m_K, m_L\}$.
(The notation $K\sqcup L=[n]$ means that $K\cup L=[n]$ and $K\cap L=\emptyset$.)
 \begin{theorem}\label{preciseCIbound}
     We have 
    $$
    \mathrm{gt}_{\m}(I)=1+\mathrm{max}\{m_{K, L}: K\sqcup L= [n]\;\text{and}\;|K|,|L|\geq 2 \}.
     $$
 \end{theorem}

% \begin{lemma}\label{notgolod}
%  If $K\sqcup L$ is a partition of $\{1\, \ldots, n\}$ with $|K|, |L|\ge 2$, then $\m^{m_{K, L}}I$ is not Golod.
% \end{lemma}
% \begin{proof}
% We have $$\prod_{i\in K} x_i^{a_i-1} \in (\m^{m_K}I):(x_i \, | \, i\in K)\subseteq (\m^{m_{K,L}}I):(x_i \, | \, i\in K)$$ and 
% $$\prod_{i\in L} x_i^{a_i-1}\in (\m^{m_L}I):(x_i\, |\, i \in L)\subseteq (\m^{m_{K,L}}I):(x_i \, | \, i\in L),$$  but $\prod_{i=1}^n x_i^{a_i-1}\notin \m^{m_{K,L}}I$. 
% \end{proof}
% This is the same argument used by Stefani in \cite{de2016products}, to show that the product of monomial ideals in 4 variables may not be Golod.

 We will need the following description of a basis of $H_{\geq 1}(\mathbf{K}^R(\mathbf{x}))$.  

\begin{proposition}\label{basischoice}
     There is a basis of $H_{\geq 1}(\mathbf{K}^R(\mathbf{x}))$ consisting elements of the form  $[udx_F]$ where $F\subseteq [n]$ is nonempty and $u\in m^{k-1}I$
     or $u=\prod_{i\in F} x_i^{a_i-1}$.        
\end{proposition}
\begin{proof}
    % Let $\Theta$ be the $k$-linear subspace of $H_{\geq 1}(\mathbf{K}^R(\mathbf{x}))$ generated by the elements whose representatives are of the above form. 
    
    We know from \cref{monomialbasisforCI}  that $H_{\geq 1}(\mathbf{K}^R(\mathbf{x}))$ has a basis consisting of elements of the form $[udx_{F}]$ where $u$ is a monomial and $F\subseteq [n]$ is nonempty. It remains to show that we can choose the monomials $u$ as in the statement. Since $udx_F$ is a cycle in the Koszul complex, we have $ux_i\in \m^kI$ for all $i\in F$. 
    This implies that $u\in \m^{k-1}I$ or $u$ is a multiple of $\prod_{i\in F}x_i^{a_i-1}$.
    Indeed, for each $i$ there exists $q_i$ such that $ux_i=hx_{q_i}^{a_i}$ with $h\in \m^k$. If $q_i\ne i$ then $x_i$ divides $h$ and it follows that $u\in \m^{k-1}I$. If $q_i=i$ for all $i\in F$, it follows that $u$ is a multiple of $\prod_{i\in F}x_i^{a_i-1}.$
    Suppose $u=M\prod_{i\in F} x_i^{a_i-1}$ with $M$ a monomial. If $M$ is a unit, our element has the required form.  Otherwise, pick $\ell \in [n]$ such that $x_{\ell}$ divides $M$. If $\ell \in F$, then $u=x_{\ell}^{a_\ell} M'\prod_{i\ne \ell, i\in F} x_i^{a_i-1}$ where $M'=M/x_{\ell}$, and it follows that $u\in \m^{k-1}I$. If $\ell \notin F$, 
    consider $y:=(u/x_{\ell})dx_{F\cup \{\ell\}}$. We have 
    \[
    \partial(y)=\pm udx_F+\sum_{i\in F}\pm (x_iu/x_{\ell})dx_{F\cup \{\ell\}\setminus \{i\}} 
    \]
    Note that $x_iu/x_{\ell} \in \m^{k-1}I$ for all $i\in F$, and therefore $[udx_F]$ is a linear combination of elements of the form $[u'dx_{F'}]$ with $u'\in \m^{k-1}I$ and $F'\subseteq [n]$.
    
    % Here, $ x_iu/x_l \in \mathfrak m^{k-1}I$. Therefore, we have \[
    % [udx_F]=\sum_{i\in F}\pm [(x_iu/x_l)dx_{F\cup \{l\}\setminus \{i\}}]\in \Theta
    % \]

    % to prove the above claim it is enough to show that an element in $H_{\geq 1}(\mathbf{K}^R(\mathbf{x}))$ represented by a monomial cycle is in the subspace $\Theta$. Let $udx_F$ be a monomial cycle in the Koszul complex and assume that $u\notin \mathfrak m^{k-1}I$. Now, as $udx_F$ is a cycle we have $u\in (\mathfrak m^k I:(x_i|i\in F))$. So, for each $i\in F$, there exists $q_i\in [n]$ such that \[
    % x_{q_i}^{a_{q_i}}|x_iu\quad \text{and}\quad \deg(x_iu/x_{q_i}^{a_{q_i}})\geq k 
    % \]
    % Now, if $q_i\neq i$, then \[\deg(x_iu/x_{q_i}^{a_{q_i}})\geq k\implies \deg(u/x_{q_i}^{a_{q_i}})\geq k-1\implies u\in \mathfrak m^{k-1}I\]
    % which is a contradiction. Therefore, $x_i^{a_i-1}|u$ for all $i\in F$. Now, if $x_i^{a_i}|u$ for some choice of $i\in F$, then $u\in \mathfrak m^{k-1}I$, which is again a contradiction. 
    % Therefore \[
    % u=M\prod_{i\in F}x_i^{a_i-1}\quad \text{where}\quad x_i\not{|}M\;\;\forall i\in F.
    % \]
    % If $M$ is a unit, then we are done, 
     
\end{proof}

\begin{lemma}\label{deg cond}
  Let $F\subseteq [n]$. The element $\prod_{i\in F} x_i^{a_i-1}dx_F$ is a cycle in $\mathbb K^R({\bf x})$ if and only if $m_F\ge k$.  
\end{lemma}
\begin{proof}
$\prod_{i\in F} x_i^{a_i-1}dx_F$ is a cycle in $\mathbb K^R({\bf x})$ if and only if for all $i\in F$ we have  \[\displaystyle x_i^a\prod_{j\in F\, \backslash \{i\}} x_j^{a_j-1} \in \m^kI,  \ \text{which is equivalent to} \prod_{ j\in F\, \backslash \, \{i\}} x_j^{a_j-1} \in \m^k\] 
This is clearly equivalent to $m_F\ge k$.
\end{proof}

We are now ready to prove \cref{preciseCIbound}
\begin{proof} 
Let $k>\mathrm{max}\{m_{K, L}: K\sqcup L= [n]\;\text{and}\;|K|,|L|\geq 2 \}$. We will use \cref{zeroproduct} to prove that $\m^kI$ is Golod.
More precisely, let $\mathbf{B}$ be a basis of $H_{\geq 1}(\mathbf{K}^R(\mathbf{x}))$ as in \cref{basischoice}. 
We claim that for every two $[z]=[udx_F]$ and $[z']=[u'dx_K]$ in $\mathbf B$, we have 
$z\cdot z'=0$

From \cref{basischoice}, we have $u\in \m^{k-1}I$ or $u=\prod_{i\in F} x_i^{a_i-1}$, and   $u'\in \m^{k-1}I$ or $u'=\prod_{i\in K} x_i^{a_i-1}$.
If $u\in \m^{k-1}I$ or $u'\in \m^{k-1}I$, then $uu'\in \m^kI$ and therefore $z\cdot z'=0$. 

Assume $u=\prod_{i\in F} x_i^{a_i-1}$ and $u'=\prod_{i\in K} x_i^{a_i-1}$. Note that we must have $|F|, |K|\ge 2$. 
Assume $F\cap K=\emptyset$.
From \cref{deg cond} we have $m_F\ge k$ and $m_K\ge k$.
Let $L:=[n]\, \backslash\,  F$. Since $F\subseteq L$, we have $m_L\ge m_F$, and therefore $m_L\ge k$.
It follows that $m_{K, L} \ge k$. This contradicts the choice of $k$, therefore $F\cap K$ must be nonempty. This implies $z\cdot z'=0$.

Now let $k\le \mathrm{max}\{m_{K, L}: K\sqcup L= [n]\;\text{and}\;|K|,|L|\geq 2 \}$. 
Let $K\sqcup L =[n]$ be a partition with $|K|,\;|L|\geq 2$ such that $k\leq m_{K,L}$. \cref{deg cond} shows that $ \left[\prod_{i\in K} x_i^{a_i-1}dx_{K}\right]$ and $\left[\prod_{i\in L} x_i^{a_i-1}dx_{L}\right]$ are elements of $H_{\geq 1}(\mathbf{K}^R(\mathbf{x})) $. Note that their product $\left[\prod_{i=1}^n x_i^{a_i-1}dx_{[n]}\right]$ is not zero in $H_{\geq 1}(\mathbf{K}^R(\mathbf{x}))$. Therefore, $ \m^kI$ is not Golod. 
\end{proof}

\section{Products of Ideals are Golod in dimension 3.}
% \begin{proposition}
%     In $P=\kk[x_1,\cdots,x_4]$, the product of any equigenerated monomial Golod ideal and the irrelevant maximal ideal is Golod. 
% \end{proposition}
% In this section, we will show that in a regular local ring of embedding dimension 3, the product of any two proper ideals is Golod.

Let $(P,\m,\kk)$ be a regular local ring of dimension 3 and let $I\subseteq \m^2$ be a proper ideal. Suppose the product structure on the positive Koszul homology over $P/I$ is trivial. Let $\mathfrak B$ be a basis for $H_{\geq 1}(\mathbf K^{P/I}(\textbf{x}))$. For $h\in \mathfrak B$, let $z_h$ be a representative for $h$. Let $\mu:\sqcup_{q\geq 1} \mathfrak B^q\xrightarrow[]{} \mathbf{K}^{P/I}(\textbf{x})$ be a map which is defined as follows: for $(h_1,\cdots,h_q)\in \mathfrak B^q$, if $q=1$, then $\mu(h_1)=z_{h_1}$, if $q=2$, then $\mu(h_1,h_2)= z$, where $\partial(z)=\overline{z_{h_1}}\cdot z_{h_2}$, if $q\geq 3$, then $\mu(h_1\ldots,h_q)=0$. One can check that $\mu$ is a trivial Massey operation. 
Therefore, a proper ideal $I\subseteq \m^2$ is Golod if and only if the product structure on the positive Koszul homology over $P/I$ is trivial, i.e., $H_{\geq 1}(\mathbf K^{P/I}(\textbf{x}))\cdot H_{\geq 1}(\mathbf K^{P/I}(\textbf{x}))=0$; see \cite{golod1962homologies}. This fact is not true in general in higher dimensions; see \cite[Theorem 3.1]{katthan2017non}.

% Let $I$ and $J$ be a proper ideals in $P$. Observe that a sufficient condition for $IJ$ to be Golod is the vanishing of one of the following maps:
%     \begin{equation}\label{golodCriterion}
%     H_{2}(\mathbf{K}^{P/IJ}(\mathbf{x}))\xrightarrow[]{}H_{2}(\mathbf{K}^{P/I}(\mathbf{x}))\quad\text{and}\quad H_{2}(\mathbf{K}^{P/IJ}(\mathbf{x}))\xrightarrow[]{}H_{2}(\mathbf{K}^{P/J}(\mathbf{x}))  
%     \end{equation}
% which are induced from the canonical projections $P/IJ\xrightarrow[]{}P/I$ and $P/IJ\xrightarrow[]{}P/J$.
%  Indeed, the vanishing of either of the above two maps allows us to choose a basis for $H_{2}(\mathbf{K}^{P/IJ}(\mathbf{x}))$ whose coefficients can be chosen to be elements of $I/IJ$ or $J/IJ$. Since a basis for $ H_1(K^{P/IJ}(\textbf{x}))$ can always be chosen to be represented by cycles with coefficients from $I$ or $J$, we have $H_{2}(\mathbf{K}^{P/IJ}(\mathbf{x}))\cdot H_{1}(\mathbf{K}^{P/IJ}(\mathbf{x}))=0$. Now, this with \cite[3.1]{vandebogert2022products} implies that $IJ$ is Golod. 

% If the ideal $I$ and $J$ are both grade 3 Gorenstein ideals in $P$, using the data of their minimal free resolution\cite{buchsbaum1977algebra} we can show that one of the maps in \cref{golodCriterion} vanishes. In general, we do not have the data of the minimal free resolution for an arbitrary ideal, but one can bypass this by embedding the ideals $I$ and $J$ into grade 3 Gorenstein ideals.

\begin{lemma}\label{mainlemma1}
    Let $(P,\m,k)$ be a regular local ring, and let $I$ and $J$ be two proper ideals in $P$. Let $\mathfrak a$ be a grade 3 Gorenstein ideal with $I\subseteq \mathfrak a$. If \begin{equation}\label{torGor}
     \Tor_{2}^P(P/IJ,\kk)\xrightarrow[]{} \Tor_{2}^P(P/\mathfrak a,\kk)   
    \end{equation}
     
    is a non-zero map, then there exists $f\in \mathfrak a$ such that $\mathfrak a\subseteq (f)+\m J$.
\end{lemma}
\begin{proof}
    Let $(G_{\cdot},\partial^{IJ}_{\cdot})$ be a minimal free resolution for $P/IJ$ over $P$, and $(F_{\cdot},\partial^{\mathfrak a}_{\cdot})$ be a minimal free resolution for $P/\mathfrak a$ over $P$. Due to the structure theorem of Buchsbaum-Eisenbud\cite{buchsbaum1977algebra} for grade 3 Gorenstein ideals, we can choose $\partial^{\mathfrak a}_2$ to be multiplication by an odd-sized alternating matrix $\Phi$, and the submaximal Pfaffians of $\Phi$ generate $\mathfrak a$. Lift the canonical projection maps $P/IJ\xrightarrow[]{}P/ \mathfrak a$ to the chain map $\alpha_{\cdot}:G_{\cdot}\xrightarrow[]{}F_{\cdot}$. Since $I\subseteq \mathfrak a$, the generators of $IJ$ are linear combinations of the generators of $\mathfrak a$ with coefficients in $J$, and therefore $\mathrm{im}(\alpha_1)\subseteq JF_1$.

Suppose \cref{torGor} is a non-zero map. This implies that $\alpha_2$ has at least one unit entry; suppose this entry is in the $ r$-th column $C_r$. Choose an $n\times n$ invertible matrix $M$ with entries in $P$ and the first column equal to $C_r$. Set $\Phi_0:=M^{T}\Phi M$, where $M^T$ is transpose of $M$. Since $M$ is invertible, the submaximal Pfaffians of both $\Phi_0$ and $\Phi$ generate the same ideal, which is $\mathfrak a$. Let $\mathbf{e}_1$ be the first canonical generator for $F_1$ and $\mathbf{b}_r$ be the $ r$-th canonical generator of $G_2$. Applying $\Phi_0$ to $\mathbf{e}_1$, we get $
\Phi_0(\mathbf e_1)= M^T \Phi M (\mathbf{e}_1)=M^T \Phi (C_r)$, and since $ \alpha_2(\textbf{b}_r)=C_r$ and $\Phi\alpha_2=\alpha_1 \partial_2$, we have 
\[
\Phi_0(\textbf{e}_1)=M^T\Phi (C_r)= M^T\Phi \alpha_2(\mathbf{b}_r)=M^T\alpha_1\partial_2(\mathbf b_r).
\] 
Since we have $\mathrm{im} (\partial_2)\subseteq \m G_1$ and $\mathrm{im}(\alpha_1)\subseteq JF_1$, we obtain \[\Phi_0(\mathbf e_1)=M^T\alpha_1\partial_2(\mathbf b_r)\subseteq M^T\alpha_1(\m G_1)\subseteq \m JF_1\]
Since $\Phi_0$ is an alternating matrix, it follows that the first row and the first column of $\Phi_0$ have entries in $\m J$. Let $f$ be the Pfaffian of $\Phi_0$ obtained by deleting the first row and column, which is one of the minimal generators of $\mathfrak a$. All the other submaximal Pfaffians are in $\m J$, and therefore $\mathfrak a\subseteq (f)+\m J$.
\end{proof}

\begin{lemma}\label{GorensteinEmbed}
    Let $I$ be a proper ideal in a Noetherian local ring $(P,\m,\kk)$. Then for any $x\notin I$, there exists a $\m$-primary Gorenstein ideal $\mathfrak a$ such that $I\subseteq \mathfrak a$ and $x\notin \mathfrak a$. 
\end{lemma}
\begin{proof}
    By the Krull intersection theorem, we have $\bigcap_{n\geq 1} (I+\m^n)=I$. Since $x\notin I$, there exists $n\geq 1$ such that $ x\notin I+\m^n$. Since $P$ is a Noetherian ring, $I+\m^n$ can be written as a finite intersection of irreducible ideals, so we can choose an irreducible ideal $\mathfrak a$ with $I+\m^n\subseteq \mathfrak a$ and $x\notin \mathfrak a$. Since $a$ is $\m$-primary, $P/\mathfrak a$ is Artinian. Now, since the zero ideal in $P/\mathfrak a$ is irreducible, we have $ \dim_{\kk} \mathrm{Soc}(P/\mathfrak a)=1$. So, $P/ \mathfrak a$ is a Cohen-Macaulay ring of type 1 and from \cite[Theorem 3.2.10]{bruns1998cohen} we know its Gorenstein.
     % To show $\mathfrak a$ is Gorenstein it is enough to show that $\dim_{\kk} \mathrm{Soc}(P/\mathfrak a)=1$. If $ \dim_{\kk} \mathrm{Soc}(P/\mathfrak a)\geq 2$, then there exists two $\kk$-linearly independent element $c_1,\;c_2\in \mathrm{Soc}(P/\mathfrak a)$. Clearly, $ (\mathfrak a +(c_1))\cap (\mathfrak a +(c_2)) = \mathfrak a $, but this contradicts the fact that $\mathfrak a$ is irreducible, therefore $\dim_\kk \mathrm{Soc}(R/I)=1$ and $\mathfrak a$ is Gorenstein. 
\end{proof}

\begin{theorem}\label{maintheorem1}
    Let $(P,\m,\kk)$ be a regular local ring of dimension 3 and $I,\; J$ be two proper ideals in $P$. Then the product ideal $IJ$ is Golod.
\end{theorem}

\begin{proof}
     Set $R=P/IJ$, and $\textbf{K}^R=\mathbf{K}^R(\textbf{x})$. From [\cite{vandebogert2022products}, 3.1], we know that $H_1(\mathbf{K}^R)\cdot H_1(\mathbf{K}^R)=0$. It remains to prove $ H_1(\mathbf{K}^R)\cdot H_2(\mathbf{K}^R)=0$. 
     
     Suppose $H_1(\mathbf K^{R})\cdot H_2(\mathbf K^{R})\neq 0$. A basis for $ H_1(\mathbf K^{R})$ can be chosen to be represented by cycles of type $f\cdot v$, where $f\in I$ and $v\in \mathbf K^{R}_1$ is such that $\partial(v)\in J$. Therefore, we may choose $[z_1]=[fv]$ with $f\in I$ and $v\in \mathbf K^{R}_1$, and $[z_2]\in H_2(\mathbf K^{R})$ such that $[z_1\cdot z_2]\ne 0$. Similarly, a basis for $H_1(\mathbf K^{R})$ can be chosen to be represented by cycles of the form $gu$ with $g\in J$ and $u\in \mathbf K^{R}_1$ such that $\partial (u)\in I$. Therefore, we can also pick $[z_1']=[gu]$ with $g\in J$ and $u\in \mathbf K^{R}_1$ such that $[z_1'z_2]\ne 0$.

     % pick $u,v\in \mathbf K^{R}_1$ and a cycle $z\in  \mathbf K^{R}_2$ such that $\partial(u)=f\in I$, $\partial(v)=g\in J$, and $[gu]\cdot [z]=[fv]\cdot [z]\neq 0$. 
     
     Since $K^{R}_3$ is a free module of rank one, we have  $[v\cdot z_2]=[rdx_{123}]$ and $[u\cdot z_2]=[sdx_{123}]$ for some $r, s\in R$.
     We claim that $r\notin J$ and $s\notin I$.
     Indeed, if $r\in J$, then $ [fv]\cdot [z_2]=[frdx_{123}]=0$, and if $s\in I$, then $ [gu]\cdot [z_2]=[gsdx_{123}]=0$, contradicting our assumption that $[z_1z_2]\ne 0$ and $[z_1'z_2]\ne 0$.
     
      From \cref{GorensteinEmbed} we have $\m$-primary Gorenstein ideals $\mathfrak a$ and $\mathfrak b$ such that $I\subseteq \mathfrak a$ with $s\notin \mathfrak a$ and $J\subseteq \mathfrak b$ with $r\notin \mathfrak b$. The following maps induced by the canonical projections are not zero:\[
    H_{2}(\mathbf{K}^R)\xrightarrow[]{\;\;\kappa_1\;\;}H_{2}(\mathbf{K}^{P/\mathfrak a}(\textbf{x}))\quad\text{and}\quad H_2(\mathbf{K}^R)\xrightarrow[]{\;\;\kappa_2\;\;}H_{2}(\mathbf{K}^{P/\mathfrak b}(\textbf{x})).
    \] 
    More precisely, we claim that  $ \kappa_1([z_2])\ne 0$ and $\kappa_2([z_2])\ne 0$. Using $\overline{\;w\;}$ to represent the image of an element $w\in\mathbf{K}^R$ along the canonical projection $\mathbf{K}^R\xrightarrow[]{} \mathbf{K}^{P/\mathfrak a}(\textbf{x})$, suppose by way of contradiction that $\overline{[z_2]}=0$. Then there exists $\ell\in \mathbf{K}^R_3$ such that $\partial(\overline{\ell})=\overline{z_2}$ in $ \mathbf{K}^{P/\mathfrak a}(\textbf{x})$. Since $\overline{u}\overline{\ell} \in \mathbf K_4^{P/\mathfrak a}=0$, we have $\partial(\overline{u})\overline{\ell}-\overline{u}\partial(\overline{\ell})=\partial(\overline{u}\overline{\ell})=0$. Recalling that $\partial(u)\in I\subseteq \mathfrak a$, which implies $\partial(\overline{u})=0$, we now have 
    $$
\overline{u}\overline{z_2}=\overline{u}\partial(\overline{\ell}) =\partial(\overline {u})\overline{\ell}=0.
    $$
    On the other hand, $[\overline{u}\overline{z_2}]=[\overline{sdx_{123}}],$ which is not zero since $s\notin \mathfrak a$. This is a contradiction.
    
    % Due to degree reasons, we have $$\partial(\overline{u}\overline{\ell})=0\implies \partial(\overline{u})\overline{\ell}-\overline{u}\partial(\overline{\ell})=0$$ 
    % Now, since $ \partial(\overline{u})=\overline{f}\in I\subseteq \mathfrak a$, and $\overline{u}\partial(\overline{\ell})=\overline{u}\overline{z_2}=\overline{sdx_{123}}$, we have $\partial(\overline{u})\overline{\ell}-\overline{u}\partial(\overline{\ell})=- \overline{sdx_{123}}=0$, but this is a contradiction as $ s\notin \mathfrak a$, therefore $\kappa_1$ is not a zero map, and similarly $\kappa_2$ is not a zero map.
    
    Now, \cref{mainlemma1} implies that there exists $\gamma \in \mathfrak a$ and $\beta \in \mathfrak b$ with $\mathfrak a\subseteq (\gamma)+\m J $ and $\mathfrak b\subseteq (\beta)+\m I $. This implies that $\mathfrak a+\mathfrak b = (\gamma,\beta)+\m( \mathfrak a+\mathfrak b)$. It follows by Nakayama's lemma that $\mathfrak a+\mathfrak b = (\gamma,\beta)$, which is a contradiction because $\mathfrak a + \mathfrak b$ is an $\m$-primary ideal and therefore it has grade 3.
\end{proof}

\begin{center}
    \section*{Acknowledgments}
\end{center}
    Experiments with Macaulay2 \cite{M2} helped us formulate results. Some of the ideas in Section 6 grew out of a discussion with ChatGPT \cite{OpenAI2026ChatGPT}. We thank Liana \c{S}ega for helpful conversations.

\bibliographystyle{abbrv}
\bibliography{references}

\end{document}